\documentclass{amsart}
\usepackage{amssymb,latexsym,amsmath}
\newtheorem{theorem}{Theorem}[section]
\newtheorem{lemma}{Lemma}[section]
\newtheorem{proposition}{Proposition}[section]

\newtheorem{remark}{Remark}[section]
\numberwithin{equation}{section}

\begin{document}
\title{Parabolic Iterated Function Systems with Applications to the Backward Continued Fractions}
\author{Andrei E. Ghenciu}
\date{10/01/2006}
\address{Department of Mathematics and Statistics\\
University of Alaska at Fairbanks\\
Fairbanks, Alaska, 99775, U.S.A.}
\email{ffeag@uaf.edu}
     
\maketitle
\begin{abstract}
To the Renyi or backward continued fraction transformation we associate
a parabolic iterated function system whose limit set has Hausdorff dimension 1.
We show that the Texan Conjecture holds, i.e. for every $t \in [0,1]$ there exists 
a subsystem whose limit set has Hausdorff dimension $t$.
\end{abstract}  
\section{Introduction}
In \cite{MU3}, Theorem 6.3., it was shown that if $S$=$\{\varphi_i\}_{i \in I}$ 
is a conformal iterated function system  and the index set $I$ is infinite, then for every
$0 \textless t \textless \theta$, there exists $A \subset I$, so that $\dim_H(J_A)$,
 the Hausdorff dimension of $J_A$, the limit set 
of the subsystem generated by $A$ is $t$; (here $\theta$ is the finiteness parameter of $S$).
For the standard continued fraction system $\theta=1/2$ and the authors conjectured that 
there is a subset A of the positive integers so that $\dim_H(J_A)=t$, 
for every $t \in [0,1]$. Later this was called the Texan Conjecture and it was an open problem 
for several years. In \cite{KS} it was proved that the Texan Conjecture was true. On the other hand, we note that there exists conformal iterated function systems $S$ and numbers $t$ between $0$ and  
the Hausdorff dimension of $J_S$, so that no subsystem of $S$ has Hausdorff dimension $t$. 
Such an example consisting of similarities was constructed in \cite{MU3}.  
\\
In this paper we prove the analogous Texan Conjecture for the backward or Renyi continued fraction system. In order to prove this, we must first construct the iterated function system corresponding 
to the backward continued fractions. This system turns out to be a parabolic iterated 
function system and this makes the analysis somewhat more involved. Moreover, we prove for a general conformal iterated function system that
if a subsystem has Hausdorff dimension $t$, then there is a regular subsystem with dimension $t$. 
The theory of general parabolic iterated function systems was developed in \cite{MU1}. 
In their work, the authors associate with a parabolic system S an (always infinite) hyperbolic conformal system S*,
whose limit set $J_{S*}$ may differ from $J_S$ by at most a countable set. This associated hyperbolic system 
is our main tool in studying the parabolic backward continued fraction system. Throughout this
paper we will make use of the notations and settings from \cite{MU1}. 

\section{Preliminaries}

In this section we collect some definitions and results from \cite{MU1} which are used throughout
this paper.  
\\
Let $X$ be a compact connected subset of a Euclidean space $\mathbf{R}^d$. Suppose that we have countably 
many conformal maps $\varphi_n:X \to X$, $n \in I$, where $I$ has at least two elements 
satisfying the following conditions:
\\
   (1)  (Open set condition) $\varphi_n(Int(X))\bigcap\varphi_m(Int(X))= \O$ $\forall m \ne n$
\\   
   (2)    $\left\|\varphi_{i} ' \right\|\textless 1$ everywhere, except for finitely many pairs 
($i$,$x_i$),
$i \in{I}$,
for which $x_i$ is the unique fixed point of $\varphi_i$ and 
$|\varphi_i(x_i)|=1$. Such pairs and indices will be called parabolic. All other
indices will be called hyperbolic.
\\
(3)  $\forall n \geq 1, \forall \omega=(\omega_1,...,\omega_n) \in {I^n}$, 
if $\omega_n$ is a hyperbolic index or $\omega_{n-1} \ne \omega_n$, then $\varphi_{\omega}$
extends conformally to an open connected subset $V \subset \mathbf{R}^d$ and maps $V$ into itself.    
\\
(4) If $i$ is a parabolic index, then $\bigcap_{n \geq 0} \varphi_{i^n}(X)=\{x_i\}$ 
and the diameters
of the sets $\varphi_{i^n}(X)$ converge to zero.
\\
(5) (Bounded Distortion Property) 
$\exists K \geq 1  \forall n \geq 1 \forall \omega=(\omega_1,...,\omega_n)
\in {I^n}  \forall x,y \in V$ if $\omega_{n-1} \ne \omega_n$ or $\omega_n$ 
is a hyperbolic index,
then 
\[
|\varphi_{\omega}'(y)| \leq K |\varphi_{\omega}'(x)|  
\]
(6)   $\exists s\textless 1  \forall n\geq 1  \forall \omega \in{I^n}$, if $\omega_n$
is a hyperbolic index or if $\omega_{n-1} \ne \omega_n$, 
then 
$$\left\| \varphi_{\omega}'\right\| \leq s$$
(7) (Cone condition) There exists $\alpha,l >0$ 
such that for every $x \in \partial X \subset \mathbf{R}^d$
there exists an open cone $Cone(x,\alpha,l) \subset Int(X)$ with vertex $x$, central angle 
of Lebesgue measure $\alpha$ and altitude $l$. 
\\
(8) There are two constants $L \geq 1$ and $\alpha \geq 0$ so that
\[
||\varphi_i'(y)|-|\varphi_i'(x)|| \leq L ||\varphi_i'|| |y-x|^ \alpha
\] 
for every $i \in I$ and
for every $x,y \in V$. 
\\
Let $\Omega$ be the set of all parabolic indices. 
If $\Omega \ne \O$, we will call the system
$\{ \phi_n:n \geq1 \}$ parabolic. The elements of $I \setminus \Omega$ are called hyperbolic. 
By $I^*$ we denote
the set of all finite words with alphabet $I$ and by $I^\infty$ we denote all the infinite 
sequences with terms in $I$.
\\
The limit set 
$J=J_S$ of the system 
$S$=$\{\varphi_i\}_{i \in I}$ is 
$J=\pi(I^\infty)$, where 
$\pi:I^\infty \to X$ is defined by 
$\{\pi(\omega)\}=\bigcap_{n \geq 0} \Phi_{\omega|n}(X)$.  
\\
For the definition of the topological pressure, finiteness parameter, Hausdorff dimension of the limit
set and for the main results concerning them for iterated function systems ,
see \cite{MU1} and \cite{MU2}.
\\
Next we briefly describe the associated hyperbolic system for a parabolic system.
Consider the system $S^*$ 
generated by the maps of the form $\varphi_{i^nj}$, where 
$n \geq1$, $i \in \Omega$, $i \ne j$ and by the maps $\varphi_k$, $k \in I \setminus \Omega$.
\\
Here are the main results we need in this paper from \cite{MU1}.
\\
Theorem: The system $S^*$ is a hyperbolic iterated function system.
\\
Theorem: The limit set $J$ and $J^*$ of the system $S$ and $S^*$ respectively differ 
only by a countable set; $J^* \subset J$ and $J \setminus J^*$ is countable.
\\
Theorem: The parabolic system $S$ is regular iff the associated system $S^*$ is regular.         

\section{Results  for General Conformal Iterated Function Systems}
Our first two theorems will be used throughout the paper. 
For a Conformal Iterated Function System $S$=$\{\varphi_i\}_{i \in I}$, 
and $A\subset I$, $\theta_A$ is the finiteness parameter of the subsystem generated by $A$,
$P_A(t)$ is the topological pressure and
$\lambda_A(t)= e^{P_A(t)} = \lim_{n \to \infty}e^{(\ln \Phi_{A,n}(t))/n}$, where
$\Phi_{A,n}(t)=\sum_{\omega \in A^n}\left\|\varphi_\omega ' \right\|^t$.
\\
The next two results are inspired by \cite{MU3}, Theorem 6.3. and \cite{KS}, Lemma 2.1., and 
allow us to estimate $\lambda_{A \cup \{b\}}(t)$ by $\lambda_A(t)$ for $A$ a subset of $I$ and $t$ a real number.   
\begin{theorem}
Let $S$=$\{\varphi_i\}_{i \in I}$, be a regular conformal iterated function system,
with $I \subset\mathbf{N}$. Let $A \subset I$ and let $b \in I$.
Let $p_b$ a positive real number so that
\begin{equation}
    \left\|\varphi_{\omega b \varpi} ' \right\|  \leq p_b \left\|\varphi_{\omega\varpi}  ' \right\|   
\end{equation}
for any words $\omega$ and $\varpi$. Then
\begin{equation}
\lambda_{A \cup \{b\}}(t) \leq \lambda_A(t) + p_b ^t
\end{equation}
for every $t \in [0,d]$.
 
\begin{proof}
Let $t \in [0,d]$ be given.
We have:
\[
\lambda_A(t) = e^{P_A(t)} = \lim_{n \to \infty}e^{(\ln \Phi_{A,n}(t))/n} = 
\lim_{n \to \infty}(\Phi_{A,n}(t))^{1/n} = \lim_{n \to \infty}
(\sum_{\omega \in A^n}\left\|\varphi_\omega ' \right\|^t)^{1/n}.
\]
Let $\tilde{A}=A \cup \{b\}$. Thus,
\[
  \lambda_{\tilde{A}}(t)=\lim_{n \to \infty}(\sum_{\omega \in \tilde{A}^n}\left\|\varphi_\omega ' \right\|^t)^{1/n}.
\]
Let $\epsilon > 0$. We show that:
\[
   \lambda_{\tilde{A}}(t) \leq \lambda_A(t) e^{\epsilon} + p_b ^t.
\] 
To this end choose $N > 0$ so that for any $n \geq N$ 
\[
   \ - \epsilon < \frac{\ln(\Phi_{A,n}(t))}{n} - P_A(t) < \epsilon.
\]
Manipulating further, we have for $n \geq N$,   
\[   
   \ e^{-\epsilon  n} < \frac{\Phi_{A,n}(t)}{(\lambda_A(t) )^n} < e^{\epsilon  n}.
\]
There exists $B > 0$, so that for all positive integers n:
\begin{equation}
\ B^{-1} e^{-\epsilon  n} < \frac{\Phi_{A,n}(t)}{(\lambda_A(t) )^n} < B e^{\epsilon  n}.
\end{equation}
Now, using (3.1) and (3.3) we get: 
\[
 \lambda_{\tilde{A}}(t)=\lim_{n \to \infty}
(\sum_{\omega \in \tilde{A}^n}\left\|\varphi_\omega ' \right\|^t)^{1/n} 
\]
\[
 \leq \lim_{n \to \infty} 
(\sum_{j=0}^{n} \binom{n}{k} p_b^{tj} \lambda_A(t)^{n-j} \frac{\Phi_{A,n-j}(t)}{(\lambda_A(t))^{n-j}})^{1/n} 
\]
\[
 \leq \lim_{n \to \infty}
(B\sum_{j=0}^{n} \binom{n}{k} p_b^{tj} \lambda_A(t)^{n-j} e^{\epsilon(n-j)} )^{1/n} 
\]
\[
=\lim_{n \to \infty} B^{1/n} ( e^\epsilon \lambda_A(t) + p_b^t) = e^\epsilon \lambda_A(t) + p_b^t
\]
Letting $\epsilon$ go to 0, we get (3.2) and this ends the proof.
\end{proof} 
\end{theorem}
Using the same ideas as in the proof of the previous theorem,
one can prove the following:
\begin{theorem}
Let $S$=$\{\varphi_i\}_{i \in I}$ be a regular conformal iterated function system, where $I \subset \mathbf{N}$.
Let $A \subset I$ and let $b \in I$. Let $r_b$ be a positive real number so that
\[
    \left\|\varphi_{\omega b \varpi} ' \right\|
    \geq r_b \left\|\varphi_{\omega\varpi}  ' \right\|   
\]
for any words $\omega$ and $\varpi$. Then
\[
   \lambda_{A \cup \{b\}}(t) \geq \lambda_A(t) + r_b ^t,
\] 
for every $t \in [0,d]$.
\end{theorem}
\begin{remark}
Under the hypothesis of Theorem 3.1., the existence of a $p_b$ for every $b$ is 
guaranteed by the Bounded Distortion Property. In particular, $p_b$ can be taken to be $K^2\left\|\varphi_b '\right\|$.
Similarly, under the hypothesis of theorem 3.2., $r_b$ can be taken to be $K^{-2}\left\|\varphi_b '\right\|$.
\end{remark}
If $S$=$\{\varphi_i\}_{i \in I}$ is a conformal iterated function system, then by the Hausdorff Dimension
spectrum of $S$, or $spec_{HD}S$, we mean the set of all $t$, so that t is the Hausdorff Dimension of 
a limit set generated by a subset of $I$. 
\\
The next result improves Theorem 6.3. in \cite{MU3}.  
\begin{proposition}
Every $t$ in the Hausdorff Dimension spectrum of a conformal iterated function system is the Hausdorff Dimension of 
a \emph{regular} subsystem. 
\begin{proof}
Let $S$=$\{\varphi_i\}_{i \in I}$ be our conformal iterated function system. 
Without loss of generality we may assume $I=\mathbf{N}$. Observe first that for every
set $A\subset \mathbf{N}$ so that $\mathbf{N} \setminus A$ is infinite and for every $\epsilon \textgreater 0$
there exists $i\in {\mathbf{N} \setminus A}$ so that $dim_H(J_{A\cup\{j\}})\leq dim_H(J_A) + \epsilon$, for every $j\geq i$.
This follows from Theorem 3.1.
\\ 
Let us fix $0\leq t \leq \theta$. We show that there exists a regular subsystem with its limit set of Hausdorff dimension $t$.
Let $I_1=\{1\}$ and suppose $I_n$ is constructed and $dim_H(J_{I_n}) \textless t$.
By the previous claim, there exists $k>max\{I_n\}$ such that 
$dim_H(J_{I_n\cup\{k\}})<t$. Choose $k_{n+1}$ be such minimal $k$ and let $I_{n+1}=I_n\cup\{k_{n+1}\}$.
Let $A=\bigcup_{n\geq1}I_n$. 
We have that $dim_H(J_A)\leq t$.
Also, $\mathbf{N}\setminus A$ is infinite;
(if not, $dim_H(J_A)\geq \theta$ and this would be a contradiction or 
would finish the proof if $dim_H(J_A)=t=\theta$).
Now, if $dim_H(J_A)=t$ we are done. Otherwise, by the previous claim, we can choose $p\in \mathbf{N}\setminus A$ 
so that $dim_H(J_{A\cup\{p\}})<t$. There must be an $n$ so that $k_n < p < k_{n+1}$. So 
$dim_H(J_{I_n\cup\{p\}})<t$ and this contradicts the choice of $k_{n+1}$.
\\
Next, we prove that the system generated by $A$ is regular.
For this, consider $\lambda_{A}(t)=e^{P_A(t)}$. 
We claim that $\lambda_{A}(t)\geq1$. 
By way of contradiction, suppose that $\lambda_{A}(t)=c<1$. 
So $\lambda_{I_n}(t)\leq c$ for every $n\geq1$. Using Theorem 3.1.,
there must be an $n$, so that $k_{n+1}\ne k_n+1$ and if $z>k_n$, then $\lambda_{I_n\cup\{z\}}(t)<1$ 
which contradicts the choice of $k_{n+1}$. 
\\
Suppose now that $t>\theta$. If $A$ is a subsystem of $I$ with its limit set of Hausdorff 
dimension $t$, then $t>\theta_A$ and so the subsystem generated by $A$ is strongly regular, and therefore regular.
This finishes the proof.
\end{proof}
\end{proposition}
We say that a conformal iterated function system $S$=$\{\varphi_i\}_{i \in I}$
has full $HD$ spectrum if for every $t \leq dim_H(J_I)$ 
there exists $A \subset I$ so that $dim_H(J_A)=t$.
\begin{proposition}
For every $t$ in the Hausdorff Dimension spectrum of a regular conformal 
iterated function system $S$, there exists a subsystem $A$ with the Hausdorff Dimension of its 
limit set arbitrarily close to $t$ and so that the subsystem generated by $A$ 
doesn't have full HD spectrum.
\begin{proof}
Let $S$=$\{\varphi_i\}_{i \in \mathbf{N}}$ our regular conformal iterated function system
and $A \subset \mathbf{N}$ so that $dim_H(J_A)=t$. 
Without loss of generality, we can choose $A$ to be infinite. 
So let $A=\{n_1,n_2,...,n_k,...\}$ and choose $k_0$ big enough so that $t_0=dim_H(J_{A_0})$ is arbitrarily close to $t$.
\\
Choose $k_1$ so that $max \{dim_H(J_{B\cup \{n_{k_1},n_{k_1+1},n_{k_1+3},...\}})\}=t_1 < t_o$ where 
$B \subset \{n_1,n_2,...,n_{k_0-1}\}$.
Let $C=\{n_1,...,n_{k_0},n_{k_1},n_{k_1+1},...\}$.
Consider the subsystem generated by $C$. The Hausdorff Dimension of its limit set arbitrarily close to $t$
and no subsystems of $C$ has Hausdorff dimension between $t_1$ and $t_0$.
\end{proof}
\end{proposition}
Next result is inspired from a result in the beginning of the proof of Theorem 6.3 in \cite{MU3}.  
Namely, if $S$=$\{\varphi_i\}_{i \in \mathbf{N}}$, with  is a conformal
iterated function system and $A$ is a subset of $\mathbf{N}$ so that $B=\mathbf{N}\setminus A$ is infinite. then for every 
$\epsilon>0$ there exists $k\in B$ so that $dim_H(J_{A\cup\{k\}})<dim_H(J_A) + \epsilon$.  
\begin{proposition}
Let $S$=$\{\varphi_i\}_{i \in I}$ be a regular parabolic iterated function system
and let $\Omega$ be the set of all parabolic indices. 
For every $i \in I$ 
let $X_i=\bigcup_{j \in I\setminus\{i\}} \varphi_j(X)$ 
and let $\beta_i= inf\{t\textgreater0 : \sum_{n\geq 1} \left\|\varphi_{i^n}'
\right\|^t_{X_i} \textless \infty \}$. 
Let $ \beta=\max \{\beta_i,i \in \Omega \}$.
If $A\subset I$ so that $I \setminus A$ is infinite and $dim_H(J_A)\textgreater \beta$, 
then for every
$\epsilon\textgreater 0$ 
$\exists b \in I \setminus A$ so that
\[
dim_H(J_{A\cup\{b\}})\textless dim_H(J_A) + \epsilon
\]
\begin{proof}
If $A$ contains no parabolic indices, then this proposition follows immediately
from the proof of Theorem 6.3 in \cite{MU3} or from the previous theorem.
So let us suppose $A$ contains the parabolic indices $i_1$,$i_2$,...,$i_k$. 
Whenever $C \subset I$, with at least one parabolic index, 
by $C^*$ we mean the associated hyperbolic system generated by $C$. 
\\
Now, whenever we add an hyperbolic index $b$ in the parabolic system $C$, we will
have the following new maps in the associated hyperbolic system:
$$\varphi_{m_n b} ; m \in \{i_1,i_2,...,i_k\} ; n \geq 0$$
Let $M=\{i_1,i_2,...,i_k\}$. Fix $t \textgreater \beta$. We have the following:
\begin{equation}
\sum_{n\geq 0,m \in M} \left\|\varphi_{m^nb}'\right\|^t_{X_m}
\leq B \left\|\varphi_b' \right\| ^t, 
\end{equation}
where $B =\{\sum_{n \geq 0,m \in M} \left\|\varphi_{m^n}' \right\|^t_{X_m}\}$.
The rest of the proof follows now from Theorem 3.1 (where $p_b$ can be replaced by $K^2\left\|\varphi_b '\right\|$)
and from the fact that the left hand side in (3.4) goes to 0 as $b$ goes to $\infty$ .
\end{proof}
\end{proposition}
\section{The Parabolic Backward Continued Fraction System}
It is well known that every irrational number $x$ in the interval [0,1) has a unique standard continued fraction expansion of the
form:
\[
 \frac{1}{\displaystyle a_1
+\frac{1}{\displaystyle a_2
+\frac{1}{a_3 + ...}}}
\]
where each of the $a_i$'s is a positive integers.
The standard continued fraction expansion is determined by the transformation:
\[
T(x)=1/x - [1/x], x \ne 0; T(0)=0.
\]
In particular, for every $n \geq 1$,:
\[
a_n=a_n(x)=[1/T^{n-1}(x)]
\]
In \cite{MU3} the authors derived many properties of subsystems of standard continued fractions
with the use of naturally associated conformal iterated function system consisting of the maps:
$$\varphi_b:[0,1] \to [0,1]; \varphi_b(x)=1/(b+x)$$    
with $b$ positive integers.
\\
We follow similar ideas based on the backward continued fraction 
expansion which we describe next.
\\
Each irrational number $x$ in the interval [0,1) has a unique backward continued fraction expansion 
of the form:
\[
 \frac{1}{\displaystyle b_1
-\frac{1}{\displaystyle b_2
-\frac{1}{b_3 - ...}}}
\]
where each of the $b_i$'s is an integer grater than 1. As with the standard continued fractions, 
there is a naturally defined transformation $T:[0,1] \to [0,1]$ which acts as the shift on 
the backward continued fractions. It is called the Renyi map and is given by:
\[
T(x)=1/(1-x) - [1/1-x]; x \ne 0; T(0)=0  
\]
In particular, for every $n\geq1$:
\[
b_n=b_n(x)=[1/(1-T^{n-1}(x))] + 1
\]
We introduce the following system generated by:
$$\varphi_b:[0,1] \to [0,1], \varphi_b(x)=1/(b-x), b \geq 2$$.
\\
Observe that $|\varphi_b '(x) | = 1/(b-x)^2$, so $\left\|\varphi_b '\right\| \textless 1  \forall b \geq 3$.
\\
Also, $|\varphi_2 '(x)| = 1/(2-x)^2$, so $|\varphi_2 '(1)|=1$ and 
$\varphi_2(1)=1$. 
\\
By induction we get that $\varphi_{2^n}'(1)=1$ for every $n \geq 0$.
This means that the system is parabolic because these maps have a neutral fixed point.
\\
Now for every $\omega \in I^*$ we can express $\varphi_{\omega}$ as a linear fractional map:
\[
\varphi_{\omega}(x)=(p_n-xp_{n-1})/(q_n-xq_{n-1}),  
\]
where $p_n=p_n(\omega)$ and $q_n=q_n(\omega)$.
\\
We have similar relations as in the case of a standard continued fraction system, namely:
\[
|\varphi_{\omega}'(x)|=1/(q_n-x q_{n-1})^2 ; \left\|\varphi_\omega'\right\|=1/(q_n-q_{n-1})^2.
\]
Also, for every $\omega \in I^*$ and for every $n\geq3$
\[
q_n(\omega)=
\omega_n q_{n-1}-q_{n-2}
\]
Using these previous relations we can state the following:
\begin{proposition}
The system generated by the maps $\varphi_b(x)=1/(b-x), b \geq 2$ where  
$\varphi_b:[0,1] \to [0,1]$ is a regular parabolic iterated function system with exactly one
parabolic index, namely 2. The Hausdorff Dimension of the system is 1.We will call this system the
BCF system. 
\end{proposition}
\begin{remark}
One can verify that the BCF system has finiteness parameter 1/2.
Also the constant $\alpha$ from Theorem 3.2 is 1/2 and $K$(the bounded distortion constant)
can be taken to be 9.
\end{remark}
\begin{remark}
The limit set of the BCF system is exactly the set of those real numbers that 
can be represented
\[
 \frac{1}{\displaystyle b_1
-\frac{1}{\displaystyle b_2
-\frac{1}{b_3 - ...}}}
\]
where $b_i$ are positive integers greater than or equal to 2. This is in fact the set of all irrational numbers in $[0,1]$.
\end{remark}
We will demonstrate that the BCF system has full spectrum.
This contrasts with our last example. We give an example of a modified continued fraction system
with full Hausdorff dimension but not with full spectrum.
\\
First, we show directly that the subsystem of the BCF system generated by the hyperbolic indices has 
full spectrum. 
\begin{proposition}
The hyperbolic iterated subsystem generated by $I=\{3,4,5,...\}$ has 
full $HD$ spectrum.
\begin{proof}
Let us consider two arbitrary words $\omega$ and $ \tilde{\omega}$ both finite and let $b \in I$
so that $b$ is strictly greater than any letters in $\omega$ and $ \tilde{\omega}$.
\\
Then for every $x \in [0,1]$ we have
\[
|\varphi_{\omega b \tilde{\omega}}'(x)|= |\varphi_{\omega}'(\varphi_b(\varphi_{\tilde{\omega}}(x)))|
|\varphi_b'(\varphi_{\tilde{\omega}}(x))|
|\varphi_{\tilde{\omega}}'(x)| 
\]
\[
\leq |\varphi_{\omega}(\varphi_{\tilde{\omega}}'(x))||\varphi_{\tilde{\omega}}'(x)| (b-1)^{-2}
\]
\[
=|\varphi_{\omega \tilde{\omega}}'(x)| (b-1)^{-2}.
\]
So,
\[
\left\|\varphi_{\omega b \tilde{\omega}}'\right\| 
\leq  \left\| \varphi_{\omega \tilde{\omega}}' \right\|(b-1)^{-2}.
\]
Now using similar ideas and the fact that the bounded distortion constant is 4, we get:
\[
\left\|\varphi_{\omega b \tilde{\omega}}'\right\| 
\geq  \left\| \varphi_{\omega \tilde{\omega}}' \right\|(2b)^{-2}. 
\]
Observe that
\[
\sum_{j \geq b+1} (2j)^{-2} \geq 1/4(b+1) 
\]
and
\[
1/4(b+1) \geq 1/(b-1)^{-2} \forall b \geq 7.
\]
This implies that for every $b \geq 7$ and for every $t \in (0,1]$ 
\[
\sum_{j \geq b+1} (2j)^{-2t} \geq 1/(b-1)^{-2t}.
\]
Thus for every $b \geq 7$ and for every $A \subset I\cap[2,b)$ 
we have:
\[
dim_H(J_{A \cup \{b \}}) \leq dim_H(J_{A \cup \{b+1,b+2,... \}}).
\] 
Therefore, (see \cite{KS}, Lemma 2.1.) the subsystem generated by $\{6,7,8,...\}$ has full $HD$ spectrum.
\\
One can check that $dim_H(J_{\{4,5,6\}}) \textless 1/2$ so we can conclude that the
subsystem generated by $\{4,5,6,...\}$ has full $HD$ spectrum.
\\
Also, using the algorithm in \cite{mU} we can prove that:
\[
dim_H(J_{\{3,4,5,6\}}) \textless 3/4 \textless dim_H(J_{\{4,5,6,...\}})
\]
Therefore the subsystem generated by $\{3,4,5,...\}$ has full $HD$ spectrum.  
\end{proof} 
\end{proposition}
\begin{lemma}
For all finite words $\omega$ and $\tilde{\omega}$ in $I^*$ and for every $n \geq 0$ we have:
\[
\left\|\varphi_{\omega 2^n b \tilde{\omega}}'\right\| 
\geq \left\|\varphi_{\omega \tilde{\omega}}'\right\|
(n+2)^{-2} (b-3/2)^{-2}
\]
for every $b \geq 5$.
\begin{proof}
Let $F:[0,1] \to [0,1]$ by $F(x)=|\varphi_{\omega 2^n b}'(x)| / |\varphi_{\omega}'(x)|$.
\\
So $F(x)=[(q_n- x q_{n-1})/(\tilde{q_{n+1}} - x \tilde{q_n})]^2$, 
\\
where $q_n=q_n(\omega)$, $\tilde{q_{n+1}}=\tilde{q_{n+1}}(\omega 2^n b)$ and $\tilde{q_n}=\tilde{q_n}(\omega 2^n)$
Observe that putting $ c = q_n - q_{n-1}$ we have:  
\[
F(1)= c^2/(\tilde{q_{n+1}} - \tilde{q_n})^2 = c^2/ ((b-1) \tilde{q_n} - \tilde{q_{n-1}})^2.
\]
Furthermore:
\[
F(1) = c^2/[(b-2)q_n + ((b-1)n - n + 1) c ]^2 = 1/[(b-2)q_n c^{-1} + ((b-1)n - n +1)]^2
\]
\[
\geq [2(b-2) + ((b-1)n - n +1)]^{-2} = [n(b-2) + 2b - 3]^{-2}.
\]
Now  $F'(x) \leq 0$ except for possibly the case when the last letter of $\omega$ is $b$. In this case using a similar technique
we can prove that $F(0) \geq [(n+2)(b-3/2)]^{-2}$.
\\
Since $F$ can have a maximum only at the end points of the interval we conclude that:
\[
F'(x) \geq [(n+2)(b-3/2)]^{-2}.
\]
Now for every $x \in [0,1]$ we have:
\[
|\varphi_{\omega 2^n b \tilde{\omega}}'(x)|
= |\varphi_{\omega 2^n b}'(\varphi_{\tilde{\omega}}(x))|
|\varphi_{\tilde{\omega}}(x)|
 \geq[(n+2)(b-3/2)]^{-2} |\varphi_{\omega}'(\varphi_{\tilde{\omega}}(x)|
\]
So:
\[
\left\|\varphi_{\omega 2^n b \tilde{\omega}}'\right\| 
\geq [(n+2)(b-3/2)]^{-2} \left\|\varphi_{\omega \tilde{\omega}}\right\|.
\]
\end{proof}
\end{lemma}
Using similar techniques we can prove the following:
\begin{lemma}
For all finite words $\omega$ and $\tilde{\omega}$ in $I^*$ and for every $n \geq 0$ we have:
\[
\left\|\varphi_{\omega 2^n b \tilde{\omega}}'\right\| 
\leq \left\|\varphi_{\omega \tilde{\omega}}'\right\|
4(n+2)^{-2} (b-1)^{-2}
\]
for every $b \geq 5$.
\end{lemma}
We are now in position to prove, with the aid of some numerical calculations, one of the main results of this paper:
We also make use of the hyperbolic subsystem associated to any parabolic subsystem.
\begin{theorem}
The BCF system has full $HD$ spectrum.
\begin{proof}
We analyze several cases. First one: $b>5$.
\\
We want to show that for every $t \geq 1/2$
\[
\sum_{n \geq 0} [2(n+2)(b-1)]^{-2t} \leq \sum_{n\geq 0,j \geq b+1} [(n+2)(b-3/2)]^{-2t}
\]
We are done if we can prove this inequality for $t=1$. So we want to prove that
\[
\sum_{n \geq 0} [2(n+2)(b-1)]^{-2} \leq \sum_{n\geq 0,j \geq b+1} [(n+2)(b-3/2)]^{-2}
\]
At this point we will be done if we show
\[
4/(b-1)^2 \leq \sum_{j \geq b+1} (b-3/2)^2
\]
and this is true for $b \geq 6$ applying the integral test.
\\
Now we want to see what happens if $b \leq 5$. For this we look at the second level of iterates
namely at the maps of the form:
\\
$\varphi_{2^r p 2^q s}$, where $r,q \geq 0$ and $p,s$ are elements of $\{3,4,....\}$.
We use Mathematica and the following formulas:
\\   
$\left\|\Phi_{2^r p 2^q s}'\right\| = A^{-2}$
where
\\
$A=(s-1)B - C$
\\
$B=p(r+1) - r + q[p(r+1)-r-(r+1)]$
\\
$B=p(r+1) - r + (q-1)[p(r+1)-r-(r+1)]$
\\
So in case $b=5$ we have
\\
$dim_H(J_{\{2,3,4,5\}}) \textless 1.995/2$ and therefore 
$dim_H(J_A) \textless 1.995/2$
for every $A \subset \{2,3,4,5\}$.
Now for $t \textless 1.995/2$ one can check using Mathematica that
\[
\sum_{n \geq 0} [2(n+2)(b-1)]^{-1.995} \leq \sum_{n\geq 0,j \geq b+1} [(n+2)(b-3/2)]^{-1.995}.
\]
Using the previous relations and Mathematica we can prove that:
\[
dim_H(J_{\{2,3,4\}}) \leq dim_H(J_{\{2,3,5,6,...\}})
\]
\[
dim_H(J_{\{2,4\}}) \leq dim_H(J_{\{2,5,6,...\}})
\]
\[
dim_H(J_{\{3,4\}}) \leq dim_H(J_{\{3,5,6,...\}})
\]
\[
dim_H(J_{\{2,3\}}) \leq dim_H(J_{\{2,4,5,6,...\}})
\]
In conclusion we proved that for every $A \subset \{2,3,4,....\}$ finite and $b$
strictly greater than any element of $A$ we have:
\[
dim_H(J_{A\cup \{b\}}) \leq dim_H(J_{A\cup \{b+1, b+2,...\}})
\]
So we can conclude that the BCF system has full $HD$ spectrum.
\end{proof}
\end{theorem}
We give an example of a modified continued fraction system of full Hausdorff dimension which 
doesn't have full spectrum:
\\
Choose $n_1$ a positive integer such that:
\begin{equation}
\sum_{j\geq n_1}1/j^{1.9}<1/3.
\end{equation}
Choose now $n_2>n_1$ such that:
\begin{equation}
1/2^{0.97} + (1/2-1/n_2)^{0.97}>1.
\end{equation}
Since $n_2>n_1$, we also have:
\begin{equation}
\sum_{j\geq n_2}1/j^{1.9}<1.
\end{equation}
Consider the system generated by the following maps:
\\
$\varphi_1:[0,1] \to [0,1], \varphi_1(x)=1/2 + x/2$,
\\
$\varphi_2:[0,1] \to [0,1], \varphi_2(x)=(1/2-1/n_2)x + 1/n_2$,
\\
$\varphi_b:[0,1] \to [0,1], \varphi_b(x)=1/(b+n_2-3+x), b\geq3$.
\\
The system has Hausdorff dimension 1.(Theorem 4.5.9. in \cite{MU4}.)
\\
Using (4.2), the Hausdorff dimension of the subsystem generated by $\{1,2\}$ is greater than 0.97.   
\\
Using (4.3), the Hausdorff dimension of the subsystem generated by $\mathbf{N}\setminus\{1\}$ 
is less than 0.95. 
Also, the Hausdorff dimension of the subsystem generated by $\mathbf{N}\setminus\{2\}$ is less than 0.95.
\\
This proves that there is no subsystem with the Hausdorff dimension of its limit set between 0.95 and 0.97.   
\[
ACKNOWLEDGMENTS
\]
The author would like to express his gratitude to his advisor, Dr. R. Daniel Mauldin, for excellent 
guidance in this project.
 
\end{document}